\documentclass[a4paper, 12pt]{extarticle}
\usepackage{setspace}
\usepackage{mathtext}
\usepackage{graphicx}
\usepackage{amscd}
\usepackage[T2A]{fontenc}
\usepackage[english]{babel}
\usepackage{latexsym,amsfonts,amssymb,amsmath,longtable,amsthm}
\usepackage[pic]{xy}
\usepackage{bbm,makeidx}
\usepackage[centerlast,small]{caption}
\usepackage{amsmath}
\usepackage[cp1251]{inputenc}
\usepackage[usenames, dvipsnames]{color}
\newtheorem{lemma}{{\scshape Lemma}}

\newtheorem{theorem}{{\scshape Theorem}}
\newtheorem{proposition}{{\scshape Proposition}}

\newtheorem{problem}{{\scshape Problem}}
\begin{document}
\title{Reidemeister spectrum of special and general linear groups over some fields contains 1}
\author{Timur Nasybullov\footnote{KU Leuven KULAK, 
Etienne Sabbelaan 53, 8500 Kortrijk, Belgium, 
 timur.nasybullov@mail.ru}~\footnote{The author is supported by the Research Foundation -- Flanders (FWO), app. 12G0317N.}}
\date{}
\maketitle
\begin{abstract}
 We prove that if $\mathbb{F}$ is an algebraically closed field of zero characteristic which has infinite transcendence degree over $\mathbb{Q}$, then there exists a field automorphism $\varphi$ of ${\rm SL}_n(\mathbb{F})$ and ${\rm GL}_n(\mathbb{F})$ such that $R(\varphi)=1$. This fact implies that ${\rm SL}_n(\mathbb{F})$ and ${\rm GL}_n(\mathbb{F})$ do not possess the $R_{\infty}$-property. However, if the transcendece degree of $\mathbb{F}$ over $\mathbb{Q}$ is finite, then ${\rm SL}_n(\mathbb{F})$ and ${\rm GL}_n(\mathbb{F})$ are known to possess the $R_{\infty}$-property \cite{Nas3}.\\

\end{abstract}
\section{Introduction}
Let $G$ be a group and $\varphi$ be an endomorphism of $G$. Elements $x$, $y$ from $G$ are called {\it $\varphi$-conjugated} if there exists an element $z\in G$ such that $x=zy\varphi(z)^{-1}$. The relation of $\varphi$-conjugation is an equivalence relation and it divides $G$ into {\it $\varphi$-conjugacy classes}. The $\varphi$-conjugacy class of the element $x$ is denoted by $[x]_{\varphi}$. The number $R(\varphi)$ of these classes is called {\it the Reidemeister number of the endomorphism $\varphi$}. The Reidemeister number is either a positive integer or infinity and we do not distinguish different infinite cardinal numbers denoting all of them by the symbol $\infty$. The subset $\{R(\varphi)~|~\varphi \text{ is an \underline{automorphism} of } G\}$ of $\mathbb{N}\cup\{\infty\}$ is called {\it the Reidemeister spectrum of $G$} and is denoted by $Spec_R(G)$. If $Spec_R(G)=\{\infty\}$, then  $G$ is said to {\it possess the $R_{\infty}$-property}.

The problem of determining groups which possess the $R_{\infty}$-property was formulated by A.~Fel'shtyn and R.~Hill \cite{FelHil}. One of the first general results in this area was obtained by A.~Fel'shtyn, G.~Levitt and M.~Lustig, they proved that non-elementary Gromov hyperbolic groups possess the $R_{\infty}$-property  \cite{Fel, LevLus}. Another extensive result was established by A.~Fel'shtyn and E.~Troitsky, they proved that every non-amenable residually finite finitely generated group possesses the $R_{\infty}$-property \cite{FelTro2}. A big list of known classes of groups which possess the $R_{\infty}$-property can be found, for example, in \cite{FelNas}. Some aspects of the $R_{\infty}$-property, namely,  relation with  nonabelian  cohomology,  relation  with  isogredience  classes  and  relation  with  representation theory can be found in \cite{FelTro}.

The author studied conditions which imply the $R_{\infty}$-property for different linear groups over rings \cite{Nas1, Nas4} and fields \cite{FelNas, Nas3}. In particular, it was proved that if $\mathbb{F}$ is an algebraically closed field of zero characteristic which has finite transcendence degree over $\mathbb{Q}$, then every reductive linear algebraic group $G$ with nontrivial quotient $G/R(G)$, where $R(G)$ denotes the solvable radical of $G$, possesses the $R_{\infty}$-property \cite{FelNas}. Also it was proved that if $\mathbb{F}$ is a field of zero characteristic which has periodic group of automorphisms, then every Chevalley group (of normal type) over the field $\mathbb{F}$ possesses the $R_{\infty}$-property \cite{Nas3}. Some fields of zero characteristic with infinite transcendence degree over $\mathbb{Q}$ have periodic groups of automorphisms (for example, the field of real numbers $\mathbb{R}$ and the field of $p$-adic numbers $Q_p$ both have trivial group of automorphisms), however, if $\mathbb{F}$ is an \underline{algebraically closed} field of zero characteristic which has infinite transcendence degree over $\mathbb{Q}$, then it always has  an automorphism of infinite order. So, the case of linear groups over algebraically closed fields of zero characteristic with infinite transcendence degree over $\mathbb{Q}$ is absolutely not studied. 

In the present paper we study twisted conjugacy classes and the Reidemeister spectrum for special linear group ${\rm SL}_n(\mathbb{F})$ and general linear group ${\rm GL}_n(\mathbb{F})$ over an algebraically closed field $\mathbb{F}$ of zero characteristic which has infinite transcendence degree over the field of rational numbers $\mathbb{Q}$. In particular, we prove that the Reidemeister spectrum of ${\rm SL}_n(\mathbb{F})$ and ${\rm GL}_n(\mathbb{F})$ contains $1$. It means that special and general linear groups over an algebraically closed field of zero characteristic with infinite transcendence degree over $\mathbb{Q}$ do not possess the $R_{\infty}$-property.  

Note that if $\mathbb{F}$ is an algebraically closed field of non-zero characteristic, then by \cite[Theorem 10.1]{Ste} for every connected linear algebraic group over $\mathbb{F}$ there exists an automorphism $\varphi$ such that $R(\varphi)=1$. This fact together with the result of the present paper and \cite{Nas1, Nas3} gives the following theorem.

~\\
\textbf{{\scshape Theorem.}} {\it Let $\mathbb{F}$ be an algebraically closed field. Then the groups ${\rm SL}_n(\mathbb{F})$ and ${\rm GL}_n(\mathbb{F})$ possess the $R_{\infty}$-property if and only if $\mathbb{F}$ has zero characteristic and the transcendence degree of $\mathbb{F}$ over $\mathbb{Q}$ is finite.}\\

The paper is organized as follows. In Section \ref{prem} we review some known simple facts about twisted conjugacy classes in different groups and recall some classical facts from Field theory (Section \ref{field}) and Model theory (Section \ref{model}). In Section \ref{rrrspec} we study the Reidemeister spectrum for special and general linear groups and prove that if $\mathbb{F}$ is an algebraically closed field of zero characteristic which has infinite transcendence degree over $\mathbb{Q}$, then there exists a field automorphism $\varphi$ of ${\rm SL}_n(\mathbb{F})$ and ${\rm GL}_n(\mathbb{F})$ such that $R(\varphi)=1$ (Theorem \ref{mth}). In Section \ref{unitel} we study properties of field automorphisms of ${\rm GL}_n(\mathbb{F})$ with finite Reidemeister number. In particular, we prove that if $\varphi$ is a field automorphism of ${\rm GL}_n(\mathbb{F})$ with $R(\varphi)<\infty$, then every matrix from ${\rm GL}_n(\mathbb{F})$ is a product of at most $3$ matrices from $[I_n]_{\varphi}\cup[I_n]_{\varphi}^{-1}$, where $I_n$ is the identity matrix (Proposition \ref{generate}). Several questions are formulated throughout the paper.

The author is grateful to Saveliy Skresanov from Novosibirsk State University for discussions on Section \ref{unitel} and for the mentioning of result \cite[Theorem 1.1]{CHLS} which simplifies the proof of Proposition \ref{generate}. Also the author is greateful to all participants of the Seminar in Algebra and Topology of KU Leuven Kulak and especially to Karel Dekimpe for valuable discussions and suggestions.

\section{Preliminaries}\label{prem}
In this section we recall some facts about twisted conjugacy classes and review some results from Field theory and Model theory which we will use in the proof of the main results.

A group $G$ is said to be \emph{divisible} if for every element $x\in G$ and every positive integer $k$ there exists an element $y\in G$ such that $y^k=x$. For example, the \textit{additive group $\mathbb{F}^+$} and \textit{the multiplecative group $\mathbb{F}^*$} of an algebraically closed field $\mathbb{F}$ are both divisible. The following proposition about abelian divisible groups is proved in \cite[Proposition 3.2]{DekGon}.
\begin{proposition}\label{dgprop}Let $\varphi$ be an endomorphism of a divisible abelian group $G$. If
$R(\varphi)$ is finite, then $R(\varphi) = 1$.
\end{proposition}
In particular, from Proposition \ref{dgprop} it follows that the Reisemeister spectrum of a divisible abelian group $G$ is a subset of $\{1,\infty\}$. The following proposition can be found in \cite[Proposition 3.1(a)]{DekGon}.
\begin{proposition}\label{dgprop2} Every divisible abelian group  does not possess the $R_{\infty}$-property.
\end{proposition}
The following result follows immediately from Proposition \ref{dgprop} and Proposition \ref{dgprop2}.
\begin{proposition}\label{triv} Let $G$ be a divisible abelian group. Then $Spec_R(G)=\{1,\infty\}$.
\end{proposition}
\noindent\textbf{Proof.} Since by Proposition \ref{dgprop2} the group $G$ does not possess the $R_{\infty}$-property, it has an automorphism  $\varphi$ with finite Reidemeister number. By Proposition \ref{dgprop} this number is equal to $1$, i.~e. 1 belongs to $Spec_R(G)$. Since every divisible group is infinite, $R(id)=\infty$ and $\infty$ belongs to $Spec_R(G)$.\hfill$\square$

~\\
The following lemma is almost obvious. A particular case of it is proved in \cite{MubSan}.
\begin{lemma}\label{qu}Let $G$ be a group, $\varphi$ be an automorphism of $G$, $H$ be a $\varphi$-admissible normal subgroup of $G$ and $\psi$ be an automorphism of $G/H$ induced by $\varphi$. Then $R(\psi)\leq R(\varphi)$.
\end{lemma}
The following equality holds for every group $G$ and its automorphism $\varphi$.
\begin{equation}\label{conn}
\left(x\varphi(x)^{-1}\right)^{y}=\left((y^{-1}x)\varphi\left(y^{-1}x\right)^{-1}\right)\left(y^{-1}\varphi(y)\right)^{-1}
\end{equation}
From this equality, in particular, follows that the subgroup $H$ of $G$ generated by the twisted conjugacy class $[e]_{\varphi}$ of the unit element is normal in $G$.
\subsection{Facts from field theory}\label{field}
If a field $L$ is a subfield of a field $F$, then we say that the field {\it $F$ is an extension of the field $L$} and denote this by $F|L$. The minimal subfield of $F$ which contains $L$ and the set of elements $X$ is denoted by $L(X)$. If for an element $x\in F$ there exists a polynomial $f(t)\neq0$ with coefficient from $L$ such that $f(x)=0$, then we say that {\it $x$ is an algebraic  element over $L$}, otherwise we say that {\it $x$ is a transcendental element over $L$}. A field extension $F|L$ is called {\it algebraic} if every element of $F$ is algebraic over $L$. If in $F$ there exists a transcendental element over $L$, then the extension $F|L$ is called {\it transcendental}.
Elements $x_1,\dots,x_k$ of the field $F$ are called {\it algebraically independent over $L$} if there is no polynomial $f(t_1,\dots,t_k)\neq0$ with coefficients from $L$ such that $f(x_1,\dots,x_k)=0$. An infinite set $X$ of elements from $F$ is called {\it algebraically independent over $L$} if every finite subset of $X$ is algebraically independent over $L$. A maximal set of algebraically independent over $L$ elements of the field $F$ is called {\it a transcendence basis of $F$ over $L$}. A cardinality of the transcendence basis of $F$ over $L$ does not depend on this basis and is called {\it a transcendence degree of the field $F$ over $L$}. The transcendence degree of $F$ over $L$ is denoted by ${\rm tr.deg}_LF$. An extension $F|L$ is said to be {\it purely transcendental} if  there exists a transcendence basis $X$ of $F$ over $L$ such that $F=L(X)$. In this case $F$ is isomorphic to the field of rational functions over the set of variables $X$ with coefficients from $L$. Every purely algebraic extension $F$ of a prime field $L$ (i.~e. $L=\mathbb{Q}$ or $L=F_p$) is completely determined by the characteristic of $L$ and the transcendence degree of $F$ over $L$. The field $F$ is said to be {\it algebraically closed} if for every polynomial $f(t)$ of non-zero degree with coefficients from $F$ there exists an element $x$ in $F$ such that $f(x)=0$. The minimal algebraically closed field which contains $F$ is called {\it the algebraic closure of $F$} and is denoted by $\overline{F}$. For every field there exists a unique (up to isomorphism) algebraic closure. The following theorem describes an arbitrary field extension as a subsequent purely transcendental extension and then algebrac extension.
\begin{theorem}\label{ext}Let $F|L$ be a field extension. Then there exists a subfield $P$ of $F$ which contains $L$ such that $P|L$ is purely transcendental and $F|P$ is algebraic.
\end{theorem}
Note that in this theorem ${\rm tr.deg}_LF={\rm tr.deg}_LP$.
If $F$ is an algebraically closed field and $L$ is a prime subfield of $F$, then by Theorem \ref{ext} there exists a subfield $P$ of $F$ such that $P|L$ is purely transcendental and $F|P$ is algebraic. Since the algebraic closure is unique, $F\cong\overline{P}$ and therefore $F$ is completely determined by its characteristic and transcendence degree over the prime subfield. The following two theorems about isomorphisms between subfields of the field of complex numbers $\mathbb{C}$ can be found in \cite{Yale}.
\begin{theorem}\label{A} Let $\varphi: F\to F^{\prime}$ be an isomorphism between two subfields of $\mathbb{C}$. If $\alpha$ is a transcendental element over $F$, then there is an isomorphism $F(\alpha)\to F^{\prime}(\beta)$ extending $\varphi$ if and only if $\beta$ is a transcendental element over $F^{\prime}$.
\end{theorem}
\begin{theorem}\label{B}Let $\varphi:F\to F^{\prime}$ be an isomorphism between two subfields of $\mathbb{C}$. Then there exists an isomorphism $\overline{F}\to\overline{F^{\prime}}$ between algebraic closures of $F$ and $F^{\prime}$ extending $\varphi$.
\end{theorem}
The if-part of Theorem \ref{A} can be obviously generalized in the following way.
\begin{theorem}\label{A'} Let $\varphi: F\to F^{\prime}$ be an isomorphism between two subfields of $\mathbb{C}$. Let the elements $\alpha_1$, $\dots$, $\alpha_k$ be algebraically independent over $F$ and $\beta_1$, $\dots$, $\beta_k$ be algebraically independent elements over $F^{\prime}$. Then there exists an isomorphism $F(\alpha_1,\dots,\alpha_k)\to F^{\prime}(\beta_1,\dots,\beta_k)$ extending $\varphi$.
\end{theorem}

\subsection{Facts from model theory}\label{model}
{\it The signature $\Sigma$} is a triple $(\mathcal{F},\mathcal{P},\rho)$, where $\mathcal{F}$ and $\mathcal{P}$ are disjoint sets not containing basic logical symbols, called, respectively, {\it the set of function symbols} and {\it the set of predicate symbols}, and $\rho:\mathcal{F}\cup \mathcal{P}\to \{0\}\cup \mathbb{N}$ is a function, called {\it arity}, which assigns a non-negative integer  to every function or predicate symbol. A function symbol $f$ is called {\it $n$-ary} if $\rho(f)=n$. A $0$-ary function symbol is called {\it a constant symbol}.

The alphabet of a first-order logic of a signature $\Sigma=(\mathcal{F},\mathcal{P},\rho)$ consists of symbols of variables (usually $x$, $y$, $z$, $\dots$), logical operations (negation $\neg$, conjunction $\wedge$, disjunction $\vee$, implication $\to$), quantifiers (existential $\exists$, universal $\forall$), functional symbols from $\mathcal{F}$ (usually $f$, $g$, $\dots$ for symbols with positive arity, and $a$, $b$, $\dots$ for $0$-ary symbols), predicate symbols from $\mathcal{P}$ (usually $p$, $q$, $\dots$), parentheses, brackets and other punctuation symbols.

{\it The set of terms} of a first-order logic of a signature $\Sigma$ is inductively defined by the following rules: any variable is a term; for a functional symbol $f\in \mathcal{F}$ with $\rho(f)=n$ and terms $t_1$, $\dots$, $t_n$ the expression $f(t_1,\dots,t_n)$ is a term.

{\it The set of formulas} of a first-order logic of a signature $\Sigma$ is inductively defined by the following rules: for a predicate symbol $p\in \mathcal{P}$ with $\rho(p)=n$ and for terms $t_1$, $\dots$, $t_n$ the expression $p(t_1,\dots,t_n)$ is a formula; for terms $t_1$, $t_2$ the expression $t_1=t_2$ is a formula; if $x$ is a variable and $\varphi$, $\psi$ are formulas, then $\varphi\wedge\psi$, $\varphi\vee\psi$, $\varphi\to\psi$, $\neg \varphi$, $\exists x~\varphi$, $\forall x~\varphi$ are formulas. A set $\mathcal{A}$ of some formulas of a first-order logic of a signature $\Sigma$ is called {\it a theory of a signature $\Sigma$}.

 Let $\Sigma=(\mathcal{F}, \mathcal{P}, \rho)$ be a signature, $M$ be a non-empty set and $\sigma$ be a function which maps every functional symbol $f$ from $\mathcal{F}$ with $\rho(f)=n$ to $n$-ary function $\sigma(f):M^n\to M$, and maps every predicate symbol $p$ from $\mathcal{P}$ to $n$-ary relation $\sigma(p)\subseteq M^n$. Denote by $\mathcal{M}$ the pair $(M,\sigma)$. Let $s$ be a function which maps any variable to some element from $M$. {\it The interpretation $[[t]]_s$ of a term $t$ in $M$ with respect to $s$} is inductively defined by the following rules: $[[x]]_s=s(x)$ if $x$ is a variable and $[[f(t_1,\dots,t_n)]]_s=\sigma(f)([[t_1]]_s,\dots,[[t_n]]_s)$ for a functional symbol $f\in \mathcal{F}$ with $\rho(f)=n$ and terms $t_1$, $\dots$, $t_n$. {\it The truth of a formula $\varphi$ in $\mathcal{M}$ with respect to $s$} (we write $\mathcal{M}\models_s\varphi$ if $\varphi$ is true in $\mathcal{M}$ with respect to $s$) is inductively defined by the following rules:
 \begin{itemize}
 \item $\mathcal{M}\models_s p(t_1,\dots,t_n)$ if an only if $\left([[t_1]]_s,\dots,[[t_n]]_s\right)\in\sigma(p)$,
  \item $\mathcal{M}\models_s t_1=t_2$ if and only if $[[t_1]]_s=[[t_2]]_s$,
   \item $\mathcal{M}\models_s\varphi\wedge\psi$ if and only if $\mathcal{M}\models_s\varphi$ and $\mathcal{M}\models_s\psi$,
   \item $\mathcal{M}\models_s\varphi\vee\psi$ if and only if $\mathcal{M}\models_s\varphi$ or $\mathcal{M}\models_s\psi$,
    \item $\mathcal{M}\models_s\varphi\to\psi$ if and only if $\mathcal{M}\models_s\varphi$ implies $\mathcal{M}\models_s\psi$,
     \item $\mathcal{M}\models_s\neg \varphi$ if and only if $\mathcal{M}\models_s\varphi$ is not true,
      \item $\mathcal{M}\models_s\exists x~\varphi$ if and only if $\mathcal{M}\models_{s^{\prime}}\varphi$ for some $s^{\prime}$ with $s^{\prime}(y)=s(y)$ for all $y\neq x$,
       \item $\mathcal{M}\models_s\forall x~\varphi$ if and only if $\mathcal{M}\models_{s^{\prime}}\varphi$ for all $s^{\prime}$ with $s^{\prime}(y)=s(y)$ for all $y\neq x$.
 \end{itemize}
We say that $\mathcal{M}=(M,\sigma)$ is a {\it model for a theory $\mathcal{A}$} if $\mathcal{M}\models_s\varphi$ for all formulas $\varphi\in \mathcal{A}$ and all functions $s$. A theory can have no models, one model or several models (even an infinite number). 
\begin{theorem}[L\"{o}wenheim–Skolem Theorem]  If a countable theory $\mathcal{A}$ of a signature $\Sigma$ has an infinite model, then for every infinite cardinal number $\kappa$ it has a model $\mathcal{M}=(M,\sigma)$ with $|M|=\kappa$.
\end{theorem}
\section{$Spec_R({\rm SL}_n(\mathbb{F}))$ and $Spec_R({\rm GL}_n(\mathbb{F}))$ contain $1$}\label{rrrspec}
In this section we study the Reidemeister spectrum for special linear group ${\rm SL}_n(\mathbb{F})$ and general linear group ${\rm GL}_n(\mathbb{F})$ and prove that if $\mathbb{F}$ is an algebraically closed field of zero characteristic with infinite transcendence degree over $\mathbb{Q}$, then the Reidemeister spectrum of ${\rm SL}_n(\mathbb{F})$ and ${\rm GL}_n(\mathbb{F})$ contains $1$. 

If $n=1$, then the statement is simple: the group ${\rm SL}_1(\mathbb{F})$ is trivial and therefore $Spec_R({\rm SL}_1(\mathbb{F}))=\{1\}$, the group ${\rm GL}_1(\mathbb{F})$ is isomorphic to the multiplecative group $\mathbb{F}^*$ of the field $\mathbb{F}$ and since $\mathbb{F}$ is an algebraically closed field, the group $\mathbb{F}^*$ is divisible, therefore by Proposition \ref{triv} we have $Spec_R({\rm GL}_1(\mathbb{F}))=Spec_R(\mathbb{F}^*)=\{1,\infty\}$. 
In the case $n\geq2$, at first, we consider a particular case of the field $\mathbb{F}$.
\begin{theorem}\label{thmain} Let $S$ be a countable set of variables, $\mathbb{F}=\overline{\mathbb{Q}(S)}$ and $G$ be either a general linear group ${\rm GL}_n(\mathbb{F})$ or a special linear group ${\rm SL}_n(\mathbb{F})$ for $n\geq2$. Then the Reidemeister spectrum of $G$ contains $1$.
\end{theorem}
\noindent\textbf{Proof.} We will prove in details only the case $G={\rm GL}_n(\mathbb{F})$, the case 
$G={\rm SL}_n(\mathbb{F})$ is similar with small mutations. Since the transcendence degree of $\mathbb{F}$ over $\mathbb{Q}$ is countable (therefore, smaller than ${\rm tr.deg}_\mathbb{Q}\mathbb{C}$), we can think about $\mathbb{F}$ as about subfield of $\mathbb{C}$ and use Theorems \ref{A}, \ref{B} and \ref{A'} for it. Since $\mathbb{F}$ is countable, the group ${\rm GL}_n(\mathbb{F})$ is also countable. Let $X_1,X_2,\dots$ be all elements of $G$.

Let us construct a chain of algebraically closed fields $F_0\subset F_1\subset F_2\subset\dots$ and its automorphisms $\varphi_0,\varphi_1,\varphi_2,\dots (\varphi_k\in{\rm Aut}~F_k)$ such that ${\rm tr.deg}_{F_k}F_{k+1}$ is finite, $\varphi_{k+1}$ extends $\varphi_k$ from $F_k$ to $F_{k+1}$ and for all $k\geq1$ there exists the set of $n^2$ elements $t_{k,ij}\in F_k~(1\leq i,j\leq n)$ such that the matrix $T_k=(t_{k,ij})$ is invertible and  $\varphi_k(T_k)=T_kX_k$, where $\varphi_k(T_k)=(\varphi_k(t_{k,ij}))$.
We will construct this sequence inductively. Let $F_0=\overline{\mathbb{Q}}$, $\varphi_0=id$ and  suppose that we already constructed fields $F_0,F_1,\dots, F_N$ and their automorphisms $\varphi_0,\varphi_1,\dots,\varphi_N$ which satisfy the required conditions.
Let $X_{N+1}=(x_{ij})$, consider the chain of fields 
$$K_{11}\subseteq K_{12}\subseteq\dots\subseteq K_{1n}\subseteq K_{21}\subseteq\dots K_{2n}\subseteq\dots\subseteq K_{nn},$$
where $K_{11}=\overline{F_N(x_{11})}$, $K_{12}=\overline{K_{11}(x_{12})}$, $\dots$, $K_{nn}=\overline{K_{nn-1}(x_{nn})}$. So, $K_{nn}$ is a minimal algebraically closed subfield of $\mathbb{F}$ which contains the field $F_N$ and all entries of the matrix $X_{N+1}$. If $x_{11}$ is a transcendental element over $F_N$, then by Theorem \ref{A} there exists an automorphism $F_N(x_{11})\to F_N(x_{11})$ extending $\varphi_N$ which fixes $x_{11}$. By Theorem \ref{B} this automorphism can be extended to the automorphism $\psi_{11}$ of $K_{11}$. If $x_{11}$ is an algebraic element over $F_N$, then $x_{11}\in F_N$ since $F_N$ is an algebraically closed field and therefore $K_{11}=F_N$. In this case denote by $\psi_{11}$ the automorphism $\varphi_N$. So,  in both cases ($x_{11}$ is either a transcendental or algebraic element over $F_N$) there exists an automorphism $\psi_{11}$ which extends $\varphi_N$ from $F_N$ to $K_{11}$. Using the same reasons there exists an automorphism $\psi_{12}$ which extends $\psi_{11}$ from $K_{11}$ to $K_{12}$, and therefore extends $\varphi_N$ from $F_N$ to $K_{12}$. Continuing this process we conclude that there exists an automorphism $\psi=\psi_{nn}:K_{nn}\to K_{nn}$ which extends $\varphi_N$ from $F_N$ to $K=K_{nn}$. Since $K$ is obtained from $F_N$ adding at most $n^2$ algebraically independent over $F_N$ elements, ${\rm tr.deg}_{F_N}K\leq n^2$.

Since ${\rm tr.deg}_{F_k}F_{k+1}$ is finite for all $k=0,\dots,N-1$ and ${\rm tr.deg}_{\mathbb{Q}}F_0=0$, then ${\rm tr.deg_{\mathbb{Q}}}F_N$ is finite and ${\rm tr.deg_{\mathbb{Q}}}K$ is finite. Since ${\rm tr.deg_{\mathbb{Q}}}\mathbb{F}$ is infinite, ${\rm tr.deg}_{K}\mathbb{F}$ is infinite, therefore in $\mathbb{F}$ there exists a set of $n^2$ elements $t_{N+1,ij} ~(1\leq i,j\leq n)$ algebraically independent over $K$. The determinant of the matrix $T_{N+1}=(t_{N+1,ij})$ is not equal to zero since $t_{N+1,ij}~(1\leq i,j\leq n)$ are algebraically independent elements over $K$ and ${\rm det}(T)$ is a polynomial over $t_{N+1,ij}~(1\leq i,j\leq n)$ with integer coefficients. Denote by $s_{ij} ~(1\leq i,j\leq n)$ the following element of $K$
$$s_{ij}=\sum_{k=1}^{n}t_{N+1,ik}x_{kj}.$$
Since $S=(s_{ij})=T_{N+1}X_{N+1}$, then $T_{N+1}=SX_{N+1}^{-1}$ and the elements $t_{N+1,ij}~(1\leq i,j\leq n)$ can be written as linear combinations of elements $s_{ij} ~(1\leq i,j\leq n)$ with coefficients from $K$. Therefore the elements $s_{ij}~(1\leq i,j\leq n)$ are algebraically independent over $K$. By Theorem~\ref{A'} there exists an isomorphism $\theta: K(t_{N+1,ij},~1\leq i,j\leq n)\to K(s_{ij},~1\leq i,j\leq n)$ which extends the automorphism $\psi$ and maps $t_{N+1,ij}$ to $s_{ij}$. Since $\theta$ is invertible, it is an automorphism of $K(t_{N+1,ij},~1\leq i,j\leq n)$.

Denote by $F_{N+1}$ the field $\overline{K(t_{N+1,ij},~1\leq i,j\leq n)}$. By Theorem \ref{B} the automorphism $\theta$ can be extended from $K(t_{N+1,ij},~1\leq i,j\leq n)$ to $F_{N+1}$. Denote this extension by $\varphi_{N+1}$. The transcendence degree of $F_{N+1}$ over $F_N$ is finite since we obtained $F_{N+1}$ from $F_N$ adding at most $2n^2$ algebraically independent over $F_N$ elements ($t_{N+1,ij}$ and probably $x_{ij}$), therefore the field $F_{N+1}$ and the automorphism $\varphi_{N+1}$ satisfy the required conditions.

Since by construction all entries  of the matrix $X_k$ belong to $F_k$, for every scalar matrix $yI_n$ there exists a number $k$ such that $y\in F_k$, therefore $\bigcup_{k}F_k=\mathbb{F}$. Denote by $\varphi$ the map $\mathbb{F}\to\mathbb{F}$ which maps the element $x\in \mathbb{F}$ to the element $\varphi_k(x)$, where $k$ is the minimal number such that $x\in F_k$. Note that for all $m>k$ we have $\varphi_m(x)=\varphi_k(x)$.
The map $\varphi$ is obviously a bijection of $\mathbb{F}$. Moreover for any two elements $x,y\in \mathbb{F}$ there exists an integer $k$ such that $x, y, x+y, xy\in F_k$, therefore $$\varphi(x+y)=\varphi_k(x+y)=\varphi_k(x)+\varphi_k(y)=\varphi(x)+\varphi(y)$$
and similarly $\varphi(xy)=\varphi(x)\varphi(y)$. Therefore $\varphi$ is an automorphism of the field $\mathbb{F}$ and it induces an automorphism of $G$ which we denote by the same letter $\varphi$. Since $T_k^{-1}\varphi(T_k)=X_k$, every $X_k$ belongs to the $\varphi$-conjugacy class of the unit element $I_n$ and therefore $R(\varphi)=1$.\hfill$\square$

The following theorem describes the general case working with an arbitrary algebraically closed field of zero characteristic.
\begin{theorem}\label{mth}Let $\mathbb{F}$ be an algebraically closed field of zero characteristic with infinite transcendence degree over $\mathbb{Q}$, and $G$ be either a general linear group ${\rm GL}_n(\mathbb{F})$ or a special linear group ${\rm SL}_n(\mathbb{F})$ with $n\geq2$. Then $Spec_R(G)$ contains $1$.
\end{theorem}
\noindent\textbf{Proof.} We will prove in details only the case $G={\rm GL}_n(\mathbb{F})$, the case
$G={\rm SL}_n(\mathbb{F})$ is similar. If the transcendence degree of $\mathbb{F}$ over $\mathbb{Q}$ is countable, then $F\cong\overline{\mathbb{Q}(S)}$, where $S$ is a countable set of variables, and the result follows from Theorem \ref{thmain}. So let the transcendence degree of $\mathbb{F}$ over $\mathbb{Q}$ be uncountable.

Let $\Sigma=(\mathcal{F},\mathcal{P},\rho)$ be a signature with $\mathcal{\mathcal{F}}=\{+,\cdot,-,^{-1},f,0,1\}$, $\mathcal{P}=\varnothing$ and $\rho(+)=\rho(\cdot)=2$, $\rho(-)=\rho(^{-1})=\rho(f)=1$, $\rho(0)=\rho(1)=0$. For terms $t_1$, $t_2$ of the signature $\Sigma$ we will write $t_1+t_2$, $t_1\cdot t_2$, $-t_1$, $t_1^{-1}$ instead of $+(t_1,t_2)$, $\cdot(t_1, t_2)$, $-(t_1)$, $^{-1}(t_1)$, respectively, and we will write $t_1\neq t_2$ instead of $\neg(t_1=t_2)$. Let $\mathcal{A}$ be the theory consisting of the following formulas.
\begin{enumerate}
\item $\forall x~\forall y~\forall z~[x+y=y+x]\wedge[x+(y+z)=(x+y)+z]$,
\item $\forall x~\forall y~\forall z~[x\cdot y=y\cdot x]\wedge[x\cdot(y\cdot z)=(x\cdot y)\cdot z]$,
\item $\forall x~\forall y~\forall z~x\cdot(y+z)=x\cdot y+x\cdot z$,
\item $\forall x~[x+0=x]\wedge[x\cdot1=x]$,
\item $\forall x~[x+(-x)=0]\wedge[(x\neq0)\to (x\cdot x^{-1}=1)]$,
\item $1\neq0$, $1+1\neq0$, $1+1+1\neq0$, $\dots$
\item $\forall y_1~\forall y_0~[(y_1\neq0)\to(\exists x~y_1\cdot x+y_0=0)]$,\\
$\forall y_2~\forall y_1~\forall y_0~[(y_2\neq0)\to(\exists x~y_2\cdot x\cdot x+y_1\cdot x+y_0=0)]$,\\
 $\dots$
\item $\forall x~\forall y~\exists z~[(f(x)=f(y))\to(x=y)]\wedge[f(z)=x]$,
\item $\forall x~\forall~y~[f(x+y)=f(x)+f(y)]\wedge[f(x\cdot y)=f(x)\cdot f(y)]$,
\item $\forall x_{11}~\forall x_{12}\dots\forall x_{1n}~\forall x_{21}\dots \forall x_{2n}\dots\forall x_{nn}~\exists y_{11}~\exists y_{12}\dots\exists y_{1n}~\exists y_{21}\dots \exists y_{2n}\dots\exists y_{nn}$\\
$[{\rm det}(x_{ij})\neq0]\to[({\rm det}(y_{ij})\neq0)\wedge\bigwedge_{i,j}(f(y_{ij})=\sum_{k=1}^ny_{ik}\cdot x_{kj})]$.
\end{enumerate}
Note that in formula 10 the symbol ${\rm det}$ is just an abbreviated form of determinant which is the polynomial with integer coefficients and can be uniquely written in terms of $\cdot$, $+$ and $-$, and  $\bigwedge_{i,j}$ means conjunction by all $1\leq i,j\leq n$.

If there exists a model $\mathcal{M}=(M,\sigma)$ for $\mathcal{A}$, then formulas 1 and 2 say that addition and multiplication are commutative and associative,  formula 3 describes distributive law between addition and multiplication, formula 4 says that $0$ and $1$ are neutral elements with respect to addition and multiplication respectively, formula 5 says that $-x$ is an opposite to $x$ element and that $x^{-1}$ is an inverse to non-zero element $x$ element. So, all together formulas 1-5 say that $(M,\sigma(+),\sigma(\cdot))$ is a field with the zero element $\sigma(0)$ and the unit element $\sigma(1)$. The countable set of formulas 6 says that the field $M$ has zero characteristic. The countable set of formulas 7 says that $M$ is an algebraically closed field. Formula 8 says that $\sigma(f)$ is a bijection on $M$ and formula 9 says that this bijection respects addition and multiplication. So, formulas 8 and 9 together say that $\sigma(f)$ is an automorphism of the field $M$. Finally, formula 10 says that for every invertible matrix $X=(x_{ij})$ with coefficients from $M$ there exists an invertible matrix $Y=(y_{ij})$ with coefficients from $M$ such that $(\sigma(f)(y_{ij}))=YX$. In other words all together formulas 1-10 say that $M$ is an algebraically closed field of zero characteristic and $\sigma(f)$ is an automorphism of $M$ which induced an automorphism $\varphi$  of ${\rm GL}_n(M)$ with $R(\varphi)=1$.

The previous paragraph is written in assumption that the theory $\mathcal{A}$ has a model, but by Theorem \ref{thmain} it has an infinite model. Since $\mathcal{A}$ is a countable theory, by L\"{o}wenheim–Skolem theorem it has a model $\mathcal{M}=(M,\sigma)$ with $|M|=\kappa=|\mathbb{F}|$, where $\mathbb{F}$ is a field from the formulation of the theorem. Since $|M|=\kappa$ is uncountable, the transcendence degree of $M$ over $\mathbb{Q}$ is infinite and is equal to $\kappa$. Since every algebraically closed field is completely determined by its characteristic and transcendence degree over prime subfield, we have $M\cong \mathbb{F}$. Therefore there exists an automorphism $\varphi$ of the field $\mathbb{F}$ which induces an automorphism of $G={\rm GL}_n(\mathbb{F})$ with Reidemeister number equals to 1.\hfill$\square$
\begin{problem}\label{prob1}{\rm Let $\mathbb{F}$ be an algebraically closed field of zero characteristic and $\psi$ be an automorphism of $\mathbb{F}$ which induces an automorphism $\varphi$ of ${\rm SL}_n(\mathbb{F})$ (or ${\rm GL}_n(\mathbb{F})$) with $R(\varphi)<\infty$. Is it possible to restrict $R(\varphi)$ in terms of $n$? Is it true that in this case $R(\varphi)=1$?}
\end{problem}
One step which we need to do in order to solve Problem \ref{prob1}  is to understand how far can be the twisted conjugacy $[I_n]_{\varphi}$ from whole group ${\rm GL}_n(\mathbb{F})$ if $R(\varphi)<\infty$. In particular if $[I_n]_{\varphi}={\rm GL}_n(\varphi)$, then $R(\varphi)=1$.
\section{Connection between $[I_n]_{\varphi}$ and ${\rm GL}_n(\mathbb{F})$}\label{unitel}
In this section we study the twisted conjugacy class $[I_n]_{\varphi}$ of the identity matrix $I_n$ in ${\rm GL}_n(\mathbb{F})$ and compare this class with whole group ${\rm GL}_n(\mathbb{F})$ in the case when $R(\varphi)<\infty$. Some connections between properties of the twisted conjugacy class of the unit element and properties of the group itself are studied in \cite{Nas2, GonNas}.

\begin{lemma}\label{dia}Let $G={\rm GL}_n(\mathbb{F})$ be a general linear group over a field $\mathbb{F}$ and $\psi$ be an automorphism of $\mathbb{F}$ which induces an automorphism $\varphi$ of $G$. If $R(\varphi)<\infty$, then the diagonal group ${\rm D}_n(\mathbb{F})$ belongs to $[I_n]_{\varphi}$.
\end{lemma}
\noindent\textbf{Proof.} Since ${\rm det}(\varphi(A))=\psi({\rm det}(A))$, a special linear group $H={\rm SL}_n(\mathbb{F})$ is a normal $\varphi$-admissible subgroup of $G$. Moreover we have the short exact sequence of groups $1\to H\to G\to \mathbb{F}^*\to 1$. The automorphism $\varphi$ of $G$ induces the automorphism $\psi$ of $G/H=\mathbb{F}^*$ (we use the same letter $\psi$ since this automorphisms coincide on $\mathbb{F}^*$). By Lemma \ref{qu} we have $R(\psi)\leq R(\varphi)<\infty$, therefore by Proposition \ref{dgprop} we have $R(\psi)=1$.

Let $X={\rm diag}(x_1,\dots,x_n)\in {\rm D}_n(\mathbb{F})$. Since $R(\psi)=1$, for every $x_i$ there exists an element $y_i$, such that $x_i=y_i\psi(y_i)^{-1}$. Denote by $Y$ the matrix ${\rm diag}(y_1,\dots,y_n)$, then $X=Y\varphi(Y)^{-1}$, i.~e. $X$ belongs to $[I_n]_{\varphi}$.\hfill$\square$

The following result is proved in \cite[Theorem 1.1]{CHLS}.
\begin{proposition}\label{ddiag} Let $A$ be a matrix from ${\rm GL}_n(\mathbb{C})$. Then there exists a diagonal matrix $D$ from ${\rm D}_n(\mathbb{C})$ such that $AD$ has $n$ different eigenvalues.
\end{proposition}
The statement of Proposition \ref{ddiag} can be easily generalized to the case of matrices over an arbitrary algebraically closed field of zero characteristic in the following way.
\begin{proposition}\label{dddiag} Let $\mathbb{F}$ be an algebraically closed field of zero characteristic and $A$ be a matrix from ${\rm GL}_n(\mathbb{F})$. Then there exists a diagonal matrix from ${\rm D}_n(\mathbb{F})$ such that $AD$ has $n$ different eigenvalues.
\end{proposition}
\textbf{Proof.} Denote by $F$ the minimal algebraically closed subfield of $\mathbb{F}$ which contains all entries of the matrix $A$. Since $F$ is obtained from $\mathbb{Q}$ by adding of at most $n^2$ algebraically independent over $\mathbb{Q}$ elements, the value ${\rm tr.deg}_{\mathbb{Q}}F$ is finite and we can think about $F$ as about a subfield of $\mathbb{C}$ and consider $A$ as a matrix over $\mathbb{C}$. 

Denote by $f: {\rm GL}_n(\mathbb{C})\to \mathbb{C}$ the function which maps a matrix $B$ from ${\rm GL}_n(\mathbb{C})$ to the discriminant of the characteristic polynomial of $B$. If we denote by $t_1,\dots,t_n$ the set of $n$ variables and by $T={\rm diag}(t_1,\dots,t_n)$, then $g(t_1,\dots,t_n):=f(AT)$ is a polynomial with coefficients from $\mathbb{C}$ over variables $t_1,\dots,t_n$. By Proposition \ref{ddiag} there exist some complex numbers $s_1,\dots,s_n$ such that $g(s_1,\dots,s_n)\neq0$. Since $g(t_1,\dots,t_n)$ is a continuous function from $\mathbb{C}^n$ to $\mathbb{C}$ and $\overline{\mathbb{Q}}$ is dense in $\mathbb{C}$ we can assume that $s_1,\dots,s_n$ are elements from $\overline{\mathbb{Q}}\subset F\subset\mathbb{F}$.
So, we proved that there exists a diagonal matrix $D={\rm diag}(s_1,\dots,s_n)$ from ${\rm D_n(\mathbb{F})}$ such that the matrix $AD$ has $n$ distinct eigenvalues.\hfill$\square$

\begin{proposition}\label{generate}Let $\psi$ be an automorphism of an algebraically closed field $\mathbb{F}$ which induces an automorphism  $\varphi$ of ${\rm GL}_n(\mathbb{F})$ with $R(\varphi)<\infty$. Then the group generated by the elements from $[I_n]_{\varphi}$ coincides with ${\rm GL}_n(\mathbb{F})$. Moreover, every matrix from ${\rm GL}_n(\mathbb{F})$ is a product of at most $3$ matrices from $[I_n]_{\varphi}\cup[I_n]_{\varphi}^{-1}$.
\end{proposition}
\noindent\textbf{Proof.} Let $A$ be an arbitrary matrix from ${\rm GL}_n(\mathbb{F})$. By Proposition \ref{dddiag} there exists a diagonal matrix $D$ such that $AD$ has $n$ distinct eigenvalues. Therefore by Jordan's theorem there exists a matrix $X$ from ${\rm GL}_n(\mathbb{F})$ such that $X^{-1}(AD)X=B$ is diagonal. Therefore $AD=XBX^{-1}$ and since by Lemma \ref{dia} the matrix $B$ belongs to $[I_n]_{\varphi}$, from equality (\ref{conn}) follows that $AD$ belongs the group generated by $[I_n]_{\varphi}$. Moreover equality (\ref{conn}) sais that $AD$ is a product of at most $2$ matrices from $[I_n]_{\varphi}\cup[I_n]_{\varphi}^{-1}$. Since by Lemma \ref{dia} the matrix $D^{-1}$ belongs to $[I_n]_{\varphi}$, the matrix $A=(AD)D^{-1}$  is a product of at most $3$ matrices from $[I_n]_{\varphi}\cup[I_n]_{\varphi}^{-1}$.\hfill$\square$

The less number $k$ we need to present every matrix from ${\rm GL}_n(\mathbb{F})$ as a product of $k$ matrices from $[I_n]_{\varphi}\cup[I_n]_{\varphi}^{-1}$, the closer the twisted conjugacy class $[I_n]_{\varphi}$ to the whole group ${\rm GL}_n(\mathbb{F})$. So, the following problem is important to study in order to understand how far can be the twisted conjugacy clas $[I_n]_{\varphi}$ from the whole group ${\rm GL}_n(\mathbb{F})$.
\begin{problem}{\rm Does there exists a number $k<3$ such that every matrix from ${\rm GL}_n(\mathbb{F})$ can be presented as a product of at most $k$ matrices from $[I_n]_{\varphi}\cup[I_n]_{\varphi}^{-1}$? If the second part of Problem \ref{prob1} is true, then such $k$ exists and is equal to $1$.}
\end{problem}

Note, that in general the twisted conjugacy class $[e]_{\varphi}$ does not have to generate the group $G$. For example, if $G$ is abelian group and $\varphi$ is an automorphism of $G$, then $[e]_{\varphi}$ is a subgroup of $G$ and $|G:[e]_{\varphi}|=R(\varphi)$. If $G={\rm GL}_n(\mathbb{F})$ and $\varphi$ is an inner automorphism of $G$, then $R(\varphi)=\infty$ and the group generated by $[I_n]_{\varphi}$ is a subgroup of the derived subgroup ${\rm GL}_n(\mathbb{F})^{\prime}={\rm SL}_n(\mathbb{F})$ which of course does not coincide with whole ${\rm GL}_n(\mathbb{F})$.

{\footnotesize
\begin{spacing}{0.5}

\end{spacing}}

\begin{thebibliography}{1}
\bibitem{Nas2}
V.~Bardakov, T.~Nasybullov, M.~Neshchadim, Twisted conjugacy classess of the unit element, Sib. Math. J., V.~54, N.~1, 2013, 10--21.
\bibitem{CHLS}
M.~Choi, Z.~Huang, C.~Li, N.~Sze, Every invertible matrix is diagonally equivalent to a matrix with distinct eigenvalues, Linear
Algebra Appl., V.~436, N.~9, 2012, 3773--3776.
\bibitem{DekGon}
 K.~Dekimpe, D.~Goncalves, The $R_{\infty}$ property for abelian groups, Topol. Methods Nonlinear Anal., V.~46, N.~2, 2015, 773--784.

\bibitem{Fel}
A.~Fel'shtyn, The Reidemeister number of any automorphism of Gromov hyperbolic group is infinite, Zap. Naucn. Semin. POMI, V.~279, 2001, 229--241.
\bibitem{FelHil}
A.~Fel'shtyn, R.~Hill, The Reidemeister zeta function with applications to Nielsen theory and a connection with Reidemeister torsion, K-Theory, V.~8, N.~4, 1994, 367--393.
\bibitem{FelNas}
A.~Fel'shtyn, T.~Nasybullov, The $R_{\infty}$ and $S_{\infty}$ properties for linear algebraic groups, J.~Group Theory, V.~19, N.~5, 2016, 901--921.
\bibitem{FelTro2}
A.~Fel'shtyn, E.~Troitsky,  Twisted conjugacy classes in residually finite groups, arXiv:math.GR/1204.3175. 
\bibitem{FelTro}
A.~Fel'shtyn, E.~Troitsky,  Aspects of the property $R_{\infty}$, J.~Group Theory, V.~18, N.~6, 2015, 1021--1034. 
\bibitem{GonNas}
D.~Goncalves, T.~Nasybullov, On groups where the twisted conjugacy class of the unit element is a subgroup, arXiv:math.GR/1705.06842. 
\bibitem{LevLus}
G.~Levitt, M.~Lustig, Most automorphisms of a hyperbolic group have very simple dynamics, Ann. Sci. Ec. Norm. Super., V.~33, 2000, 507--517.
\bibitem{MubSan}
T.~Mubeena, P.~Sankaran, Twisted conjugacy classes in abelian extensions of certain linear groups, Canad. Math. Bull., V.~57, 2014, 132--140.
\bibitem{Nas1}
T.~Nasybullov, Twisted conjugacy classes in general and special linear groups,\\
  Algebra and Logic,  V.~51, N.~3, 2012, 220--231.

\bibitem{Nas3}
T.~Nasybullov, Twisted conjugacy classes in Chevalley groups, Algebra Logic,  V.~53, N.~6, 2014, 481--501.
\bibitem{Nas4}
 T.~Nasybullov, The $R_{\infty}$-property for Chevalley groups of types $B_l$, $C_l$, $D_l$ over integral domains,  J.~Algebra, V. 446, 2016, 489--498.
\bibitem{Ste}
R.~Steinberg, Endomorphisms of Linear Algebraic Groups, Memoirs of AMS, V.~80, 1968.
\bibitem{Yale}
P.~Yale, Automorphisms of the complex numbers, Math. Mag., V.~39, N.~3, 1966, 135--141.
\end{thebibliography}
\end{document}